# Intentionally and Unintentionally.
# On Both, A and Non-A, in Neutrosophy

Feng Liu
Department of Management Science and Engineering
Shaanxi Economics and Trade Institute (South Campus)
South Cuihua Road, Xi'an, Shaanxi, 710061, P. R. China
E-mail: liufeng49@sina.com

Florentin Smarandache
Department of Mathematics
University of New Mexico, Gallup, NM 87301, USA
E-mail: smarand@unm.edu

**Abstract:** The paper presents a fresh new start on the neutrality of neutrosophy in that "both A and Non-A" as an alternative to describe Neuter-A in that we conceptualize things in both intentional and unintentional background. This unity of opposites constitutes both objective world and subjective world. The whole induction of such argument is based on the intensive study on Buddhism and Daoism including I-ching. In addition, a framework of contradiction oriented learning philosophy inspired from the Later Trigrams of King Wen in I-ching is meanwhile presented. It is shown that although A and Non-A are logically inconsistent, but they are philosophically consistent in the sense that Non-A can be the unintentionally instead of negation that leads to confusion. It is also shown that Buddhism and Daoism play an important role in neutrosophy, and should be extended in the way of neutrosophy to all sciences according to the original intention of neutrosophy.

**Keywords:** Neutrosophy, Multi-valued logic, I-ching, Buddhism, Dao, Yin-yang, Creativity, Contradiction, Illusion, Learning

**2000 MSC**: 03B60, 03B42

## 1. Objective world and subjective world

The common confusion about the objective world is: it **is** just what we see and feel. This is however very wrong. In fact, this is rather a belief than an objective reflection, and varies among different people, because none of us can prove it. In his paper "To be or not to be, A multidimensional logic approach" Carlos Gershenson [2] has generalized proofs:

- Everything is and isn't at a certain degree. (i.e., there is no absolute truth or false);
- Nothing can be proved (that it exists or doesn't) (i.e., no one can prove whether his consciousness is right);
- I believe, therefore I am (i.e., I take it true, because I believe so).

It is something, but not that figured in our mind. This is the starting point in Daoism (F. Liu [2]).

Daodejing begins with: "Dao, daoable, but not the normal dao; name, namable, but not the normal name." We can say it is dao, but it doesn't mean what we say. Whenever we mention it, it is beyond the original sense.

Daodejing mainly deals with the common problem: "What/who creates everything in the world we see and feel?" It is dao: like a mother that bears things with shape and form. But what/who is dao? It is just unimaginable, because whenever we imagine it, our imagination can never be it (we can never completely describe it: more we describe it, more wrong we are). It is also unnamable, because whenever we name it, our concept based on the name can never be it.

Daoism illustrates the origin of everything as such a form that doesn't show in any form we can perceive. This is the reason why it says, everything comes from nothingness, or this nothingness creates everything in forms in dynamic change. Whatever we can perceive is merely the created forms, rather than its genuine nature, as if we distinguish people by their outer clothes. We are too far from understanding the nature, even for the most prominent figures like Einstein.

Therefore **Name and Non-Name coexist** pertaining to an object:

Object = both Name and Non-Name.

Then what should we do **subjectively**? Very simple: both intentionally and unintentionally. Intentional conception relates to all the connotation and extension pertaining to Name, and unintentional one to Non-Name.

- **There are alternative interpretations on Non-Name: unintentionally and negatively. This is crucial in our confusion.**

This is the contradiction between creativity and implementation, as is stated below.

2. Neutrosophy

Neutrosophy is a new branch of philosophy that studies the origin, nature, and scope of neutralities, as well as their interactions with different ideational spectra.
It is the base of *neutrosophic logic*, a multiple value logic that generalizes the fuzzy logic and deals with paradoxes, contradictions, antitheses, antinomies.

**Characteristics** of this mode of thinking:
- proposes new philosophical theses, principles, laws, methods, formulas, movements;
- reveals that world is full of indeterminacy;
- interprets the uninterpretable;
- regards, from many different angles, old concepts, systems:
showing that an idea, which is true in a given referential system, may be false in another one, and vice versa;
- attempts to make peace in the war of ideas,
and to make war in the peaceful ideas;
- measures the stability of unstable systems,
and instability of stable systems.
Let's note by <A> an idea, or proposition, theory, event, concept, entity, by <Non-A> what is not <A>, and by <Anti-A> the opposite of <A>. Also, <Neut-A> means what is neither <A> nor <Anti-A>, i.e. neutrality in between the two extremes. And <A'> a version of <A>.

<Non-A> is different from <Anti-A>.

**Main Principle:**

Between an idea <A> and its opposite <Anti-A>, there is a continuum-power spectrum of neutralities <Neut-A>.

**Fundamental Thesis of Neutrosophy:**

Any idea <A> is T% true, I% indeterminate, and F% false, where T, I, F $\subset$ ]$^-$0, 1$^+$[.

**Main Laws of Neutrosophy:**

Let <α> be an attribute, and (T, I, F) $\subset$ ]$^-$0, 1$^+$[$^3$. Then:

- There is a proposition <P> and a referential system {R}, such that <P> is T% <α>, I% indeterminate or <Neut-α>, and F% <Anti-α>.
- For any proposition <P>, there is a referential system {R}, such that <P> is T% <α>, I% indeterminate or <Neut-α>, and F% <Anti-α>.
- <α> is at some degree <Anti-α>, while <Anti-α> is at some degree <α>.

## 3. Creativity and implementation

We can model our mind in the alternation of yin and yang that is universal in everything (Feng Liu):

- Yang pertains to dynamic change, and directs great beginnings of things; yin to relatively static stage, and gives those exhibited by yang to their completion.

  In the course of development and evolution of everything yang acts as the creativity (Feng Liu) that brings new beginnings to it, whereas yin implements it in forms as we perceive as temporary states. It is in this infinite parallelism things inherit modifications and adapt to changes.

- On our genuine intelligence — **creativity** (F. Liu [2])

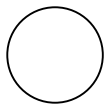

In the query about the figure on the left, whenever we hold the answer as a circle, we are inhibiting our creativity. Nor should we hold that it is a cake, a dish, a bowl, a balloon, or the moon, the sun, for we also spoil our creativity in this way. Then, what is it?

"It is nothing."

Is it correct? It is, if we do not hold on to the assumption "it is something". It is also wrong, if we persist in the doctrine "the figure is something we call nothing." This nothing has in this way become something that inhibits our creativity. How ridiculous!

**Whenever we hold the belief "it is …", we are loosing our creativity. Whenever we hold that "it is not …", we are also loosing our creativity.** Our true intelligence requires that we completely free our mind — neither stick to any extremity nor to "no sticking to any assumption or belief". This is a kind of genius or gift rather than logic rules, acquired largely after birth, e.g., through Buddhism practice. **Note that our creativity lies just between internationality and uninternationality.**

**Not (it is) and not (it is not),**
**It seems nothing, but creates everything,**
**Including our true consciousness,**
**The power of genius to understand all.**

- The further insight on contradiction compatible **learning philosophy** inspired from the Later Trigrams of King Wen of I-ching shows that:

    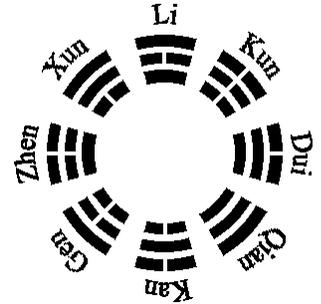

    When something (controversial) is perceived (in Zhen), it is referred (in Xun) to various knowledge models and, by assembling the fragments perceived from these models, we reach a general pattern to which fragments attach (in Li), as hypothesis, which needs to be nurtured and to grow up (Kun) in a particular environment. When the hypothesis is mature enough, it needs to be represented (in Dui) in diverse situations, and to expand and contradict with older knowledge (in Qian) to make update, renovation, reformation or even revolution in knowledge base, and in this way the new thought is verified, modified and substantialized. When the novel thought takes the principal role (dominant position) in the conflict, we should have a rest (in Kan) to avoid being trapped into depth (it would be too partial of us to persist in any kind of logic, to adapt to the outer changes). Finally the end of cycle (in Gen).

    I-ching [in Chinese: Yi Jing] means: *Yi* = change, *Jing* = scripture. It mainly deals with the creation, evolution (up and down) of everything in such perspective that everything is an outer form of a void existence, and that everything always exists in the form of unity (compensation, complementation) of opposites…

    This philosophy shows that **contradiction acts as the momentum or impetus to learning evolution. No controversy, no innovation. This is the essentially of neutrosophy** (Florentin Smarandache).

  In the cycle there is unintentionally implied throughout it:
    - Where do the reference models relating to the present default model come from? They are different objectively.
    - How can we assemble the model from different or even incoherent or inconsistent fragments?
    - If we always do it intentionally, how does the hypothesis grow on its own, as if we study something without sleep?
    - How can our absolute intention be complemented without contradiction?
    - Is it right that we always hold our intention?
        There is only one step between truth and prejudice — when the truth is overbelieved regardless of constraint in situations, it becomes prejudice.

- Is there no end for the intention? Then, how can we obtain a concept that is never finished? If there is an end, then it should be the beginning of unintentionally, as yin and yang in Taiji figure.
  Unintentionally can be alternately derived when the intention is repeated over and over so that it becomes an instinct. **It is even severe that we develop a fallacy instinct.**

**4. Completeness and incompleteness: knowledge and practice**

There always is contradiction between completeness and incompleteness of knowledge. In various papers presented by Carlos Gershenson he proves this point. Same in Daoism and Buddhism. This contradiction is shown in the following aspect:
  a) People are satisfied with their knowledge relative to a default, well-defined domain. But later on, they get fresh insight in it.
  b) They face with contradictions and new challenges in their practice and further development.

This reflects in our weakness that:
- Do we understand ourselves?
- Do we understand the universe?

What do we mean by knowledge, complete whole or incomplete? Our silliness prompts us to try the complete specifications, but where on earth are they (Gershenson [1])? Meanwhile, our effort would be nothing more than a static imitation of some dynamic process (Liu [1]), since **human understands the world through the interaction of the inter-contradictory and inter-complementary two kinds of knowledge: perceptual knowledge and rational knowledge - they can't be split apart.**

- In knowledge discovering, there are merely strictly limited condition that focus our eyes to a local domain rather than a open extension, therefore our firsthand knowledge is only relative to our default referential system, and extremely subjective possibly.
- Is it possible to reach a relatively complete piece at first? No, unless we were gods (we were objective in nature).
- Then we need to perceive (the rightness, falseness, flexibility, limitation, more realistic conception, etc.) and understand the real meaning of our previous knowledge — how: only through practice (how can we comprehend the word "apple" without tasting it?)
- Having done that, we may have less subjective minds, based on which the original version is modified, revised, and adapted as further proposals.
- Again through practice, the proposals are verified and improved.
- This cycle recurs to the infinite, in each of which our practice is extended in a more comprehensive way; the same to our knowledge.
  - Discover the truth through practice, and again through practice verify and develop the truth. Start from perceptual knowledge and actively develop it into rational knowledge; then start from rational knowledge and actively guide revolutionary practice to change both the subjective and the objective world. Practice, knowledge, again practice, and again knowledge. This form repeats itself in endless cycles, and with each cycle the content of practice and knowledge rises to a higher level.
- Through practice, we can
  - verify our knowledge
  - find the inconsistency, incompleteness of our knowledge, and face new problems, new challenges as well
  - maintain critical thought

Therefore **knowledge is based on the infinite critics and negation (partial or revolutionary) on our subjective world.** It is never too old to learn.

## 5. Conclusion

Whenever we say "it is", we refer it to both subjective and objective worlds.

We can creatively use the philosophical expression **both A and non-A** to describe both subjective world and objective worlds, and possibly the neutrality of both.

Whenever there is "it is", there is subjective world, in the sense that concepts always include subjectivity. So our problem becomes: is "it" really "it"? A real story of Chinese Tang dynasty recorded in a sutra (adapted from Yan Kuanhu Culture and Education Fund)shows that:

> Huineng arrived at a Temple in Guangzhou where a pennant was being blown by wind. Two monks who happened to see the pennant were debating what was in motion, the wind or the pennant.
>
> Huineng heard their discussion and said: "It was neither the wind nor the pennant. What actually moved were your own minds." Overhearing this conversation, the assembly (a lecture was to begin) were startled at Huineng's knowledge and outstanding views.

- When we see pennant and wind we will naturally believe we are right in our consciousness, however it is subjective. In other words, what we call "the objective world" can never absolutely be objective at all.
- Whenever we believe we are objective, this belief however is subjective too.
- In fact, all these things are merely our mental creations (called illusions in Buddhism) that in turn cheat our consciousness: There is neither pennant nor wind, but our mental creations.
- The world is made up of our subjective beliefs that in turn cheat our consciousness. This is in fact a cumulative cause-effect phenomenon.
- Everyone can extricate himself out of this maze, said Sakyamuni and all the Buddhas, Bodhisattvas around the universe, their number is as many as that of the sands in the Ganges (Limitless Life Sutra).

**References:**

C. Gershenson [1]: *Comments to Neutrosophy*, Proceedings of the First International Conference on Neutrosophy, Neutrosophic Logic, Set, Probability and Statistics, University of New Mexico, Gallup, December 1-3, 2001, http://www.gallup.unm.edu/~smarandache/CommentsToNeutrosophy.pdf